\documentclass[12pt]{amsart}
\usepackage[pdftex]{graphicx}
\newtheorem{theorem}{Theorem}[section]
\newtheorem{proposition}{Proposition}[section]
\newtheorem{corollary}{Corollary}[section]

\begin{document}

\title{Eigenvalues  of the  bilayer graphene operator with a   complex  valued  potential}

\author{ Francesco Ferrulli,  Ari Laptev  and  Oleg  Safronov }

\address{Francesco Ferrulli, Department of Mathematics, Imperial College London, SW7 2AZ, London, UK}
\email{f.ferrulli14@imperial.ac.uk}

\address{Ari Laptev, Department of Mathematics, Imperial College London, SW7 2AZ, London, UK}
\email{laptev@mittag-leffler.se}

\address{Oleg Safronov, Department of Mathematics	 and Statistics, University of North Carolina at Charlotte, Charlotte, NC 28223, USA}
\email{osafrono@uncc.edu}

\begin{abstract}
We  study the spectrum  of a system of second order differential operator $D_m$ perturbed  by a non-selfadjoint   matrix  valued  potential $V$. We  prove  that  eigenvalues  of  $D_m+V$ are located  near the edges of  the  spectrum of the unperturbed operator $D_m$.
\end{abstract}

\maketitle

\section{Statement of the main results}

Spectral properties  of  non-selfadjoint operators have been recently a subject of interest of many papers. A particular interest was  related to the location of eigenvalues of  differential operators  in the complex  plane ${\Bbb C}$. The corresponding results for  Schr\"odinger operators  can be  found in \cite{AAD}, \cite{Da}-\cite{DN}  and  in  \cite{Fr}. Some other problems  were  studied  in the papers
\cite{FrLaLiSe}-\cite{FrSi}  and \cite{LaSa}.

The operator  we  study  is  related to the  quantum  theory of a material consisting of  two layers  of  graphene.
Namely, we   consider the operator $D=D_m+V$,
where 
$$
D_m=
\begin{pmatrix}
m&4\partial_{\bar z}^2\\4\partial_{ z}^2&-m
\end{pmatrix},\qquad  \partial_{\bar z}=\frac12\Bigl(\frac{\partial}{\partial x_1}+\frac1i\frac{\partial}{\partial x_2}\Big),  \qquad \partial_{ z}=\frac12\Bigl(\frac{\partial}{\partial x_1}-\frac1i\frac{\partial}{\partial x_2}\Big),\quad m\geq 0 . 
$$
This operator  acts in the Hilbert  space $L^2({\Bbb R}^2; {\Bbb C}^2)$. The domain of $D$ is  the  Sobolev  space ${\mathcal H}^2({\Bbb R}^2; {\Bbb C}^2)$.
The potential $V$ is a not necessary self-adjoint matrix-valued function
$$
V(x)=\begin{pmatrix}
V_{1,1}(x)&V_{1,2}(x)\\V_{2,1}(x)&V_{2,2}(x),
\end{pmatrix}
$$ 
where the matrix  elements  are allowed  to take  complex  values. For the matrix $V$ we denote
$$
|V(x)| = \sqrt{\sum_{i,j=1,2} |V_{i,j}(x)|^2}.
$$
Assuming that $V$  decays  at   the infinity in some  integral sense we would like to answer the  question: "Where are  the eigenvalues of $D$  located?"

Note that since $D_m^2=\Delta^2+m^2$, the  spectrum $\sigma(D_m)$ of $D_m$  is the set $(-\infty, m]\cup [m,\infty)$. Our  results  show  that   the eigenvalues   of $D$ are located  near the edges   of the  absolutely  continuous  spectrum, i.e.  near the points  $\pm m$.
Since   the spectrum of  the unperturbed operator  has two edges,  our results  resemble some of  the  theorems of  the paper \cite{CLT} related  to the Dirac operator.  However, the main difference between the  two papers  is  that  we study a differential  operator on a plane, while  the article \cite{CLT}  deals  with  operators on a  line.

\begin{theorem}\label{t1.1} Let $k\notin \sigma(D_m)$ be an eigenvalue  of the operator $D$. Let $1<p<4/3$.
Then
$$
\frac{C_p\int_{{\mathbb R}^2} |V(x)|^p dx}{|\mu|^{p-1}}\Bigl(\sqrt{\Bigl|\frac{k-m}{k+m}\Bigr|}+\sqrt{\Bigl|\frac{k+m}{k-m}\Bigr|}+1\Bigr)^p\geq 1,\qquad \mu^2=k^2-m^2,
$$ with $C_p>0$ independent of $V$, $k$ and $m$.
In particular, if $m=0$, then
$$
|k|^{p-1}\leq 3^{p}C_p\int_{{\mathbb R}^2} |V(x)|^p dx,\qquad 1<p<4/3.
$$
\end{theorem}
\bigskip 

The next  statement tells us about what happens when  $p\to1$.

\begin{theorem}\label{t1.2} Let $k\notin \sigma(D_m)$ be an eigenvalue  of the operator $D$. Let $ \mu^2=k^2-m^2.$
Then
$$
C\Bigl(|\ln|\mu||\sup_{x\in{\mathbb R}^2}\int_{|x-y|<(2|\mu|)^{-1}} |V(y)|\,dy+\sup_{x\in{\mathbb R}^2}\int_{{\mathbb R}^2}\Bigl(1+|\ln|x-y||\Bigr)|V(y)|dy
\Bigr)+$$
$$+C\int_{{\mathbb R}^2}|V(x)|dx\Bigl(\sqrt{\Bigl|\frac{k-m}{k+m}\Bigr|}+\sqrt{\Bigl|\frac{k+m}{k-m}\Bigr|}+1\Bigr)\geq 1,
$$
where  the  constant $C>0$ is independent of $V$, $k$ and $m$.
\end{theorem}

\bigskip

Note  that   this statement   also holds  true  for $m=0.$  

\begin{corollary}Let $m=0$ and let $k\notin {\Bbb R}$ be an eigenvalue  of the operator $D$. 
Then
$$
C\Bigl(|\ln|k||\sup_{x\in{\mathbb R}^2}\int_{|x-y|<(2|k|)^{-1}} |V(y)|\,dy+\sup_{x\in{\mathbb R}^2}\int_{{\mathbb R}^2}\Bigl(1+|\ln|x-y||\Bigr)|V(y)|dy
\Bigr)+$$
$$+3C\int_{{\mathbb R}^2}|V(x)|dx\geq 1,
$$
where  the  constant $C>0$ is independent of $V$ and $k$.
\end{corollary}

In particular, we see   that  if $m=0$,  then   for  small $V$,   the eigenvalues  of  $D$  are situated in the  circle $\{k\in{ \Bbb C}:\,\,\,\,|k|<r\}$ of  radius $r$   which   has the following asymptotical    behavior 
$$
r\asymp \exp \Bigl(-\frac{C}{\int |V|dx}\Bigr),\qquad {\rm as}\qquad \int |V|dx\to 0.
$$

\medskip
\noindent
The proof of Theorems 1.1 and 1.2 are given in Section 2. 
In Section 3 we consider a special case where $V = i W^2$, $W= W^*$, 
 In this case we can get a more precise information about location of the complex eigenvalues, see Theorem 3.1. It is interesting to note that if $m=0$ (no gap in the continuous spectrum), then perturbations by such matrix-functions do not create any complex eigenvalues. Here we have similarities with the result obtained for the one dimensional Dirac operators in \cite{CLT}.

\bigskip
\section{Proofs  of the main results}

In order to prove our   main results we need  the Birman-Schwinger  principle  formulated   below.
\begin{proposition} Let $V=W_2 W_1,$  where $W_1$  and $W_2$ are  two  matrix-valued decaying  functions.
A point $k\in {\Bbb C}\setminus \sigma(D_m)$ is an eigenvalue of $D$  if and  only if  $-1$ is an eigenvalue of the operator
$$X(k):=W_1(D_m-k)^{-1}W_2.
$$
In  particular, if $k\in  {\Bbb C}\setminus  \sigma(D_m)$  is  an eigenvalue of $D$ then $||X(k)||\geq1$.
\end{proposition}

\noindent
The proof of this  statement is standard and it is left  to the reader as an exercise.  

\medskip
\noindent
Below we  always denote
$$
W=\sqrt{V^*V}
$$
 and  use the Birman-Schwinger  principle with $W_1=W$ and $W_2 =V W^{-1/2}$. 

\bigskip
\noindent
{\it Proof of Theorem~\ref{t1.1}}.  Since $$(D_m-k)^{-1}=(D_m+k)(D_m-k)^{-1}(D_m+k)^{-1}=(D_m+k)(D_m^2-k^2)^{-1},
$$
it is easy to see  that
\begin{equation}\label{resolvent}
(D_m-k)^{-1}=(m\gamma_0+k-\mu)(\Delta^2-\mu^2)^{-1}+(D_0-\mu)^{-1},
\end{equation}
where 
$$
D_0=
\begin{pmatrix}
0&4\partial_{\bar z}^2\\4\partial_{ z}^2&0
\end{pmatrix},\qquad \gamma_0=\begin{pmatrix}1&0\\0&-1
\end{pmatrix}.
$$
One can also  note that  the last term in the right hand  side of \eqref{resolvent} can be rewritten in the form
\begin{equation}\label{resolventm=0}
(D_0-\mu)^{-1}=(D_0+\mu)(\Delta^2-\mu^2)^{-1}.
\end{equation}
The operator $(\Delta^2-\mu^2)^{-1}$  is an integral operator with the  kernel
$$
g_k(x,y)=\frac{i}{8\mu}\Bigl(H(\sqrt{\mu} r)-H(i\sqrt{\mu} r)\Bigr),
$$
where $H(z)=H^{(1)}_0(z)$  is  the Hankel function of first kind  and $r=|x-y|$.   It is a simple consequence of the fact  that
$$
(\Delta^2-\mu^2)^{-1}=\frac1{2\mu}\Bigl((-\Delta-\mu)^{-1}-(-\Delta+\mu)^{-1}\Bigr).
$$
The  kernel of $(-\Delta-\mu)^{-1}$  is $4^{-1}iH(\sqrt \mu r)$.
Another   useful representation of $g_k(x,y)$ follows  from  the  fact  that   the kernel of  $(-\Delta-\mu)^{-1}$ equals  (see \cite{GS})
$$
(2\pi)^{-1}K_0(-i\sqrt{\mu} |x-y|), 
$$
where
$$
K_0(z)=\frac{e^{-z}}{\Gamma(1/2)}\sqrt{\frac{\pi}{2z}}\int_0^\infty  e^{-t}t^{-1/2}\bigl(1+\frac{t}{2z}\bigr)^{-1/2}dt,\qquad  |\arg z|<\pi.
$$
Let us  define
$$
G(z)=H(z)-H(iz).
$$
We need  to know the  behaviour  of the  function $G$ only in the region $0<\arg z<\pi/2$, where we have
$$
|G(z)|+|G'(z)|+|G''(z)|\leq \frac{C}{\sqrt{|z|}},\qquad {\rm if}\qquad |z|>1/2.
$$
The behaviour  of  the function $G$  near  zero is determined  by  the expansion of the  Hankel  function in the  neighbourhood of $z=0$.
It turns out that
$$
|G(z)|\leq C,\qquad |G'(z)|\leq C_1 |z|\ln|z|^{-1},\qquad |G''(z)|\leq C_1 \ln|z|^{-1},\qquad {\rm if}\qquad |z|<1/2.
$$
Let $\rho_\mu(|x-y|)$ be the kernel of the integral operator $(D_0-\mu)^{-1}$
\begin{equation*}
\rho_\mu(|x-y|) =  \frac{i}{8\mu}
\begin{pmatrix} 
\mu G(\sqrt{\mu}|x-y|) & \partial_{\bar{z}}^2 G(\sqrt{\mu}|x-y|)
\\\partial_{z}^2 G(\sqrt{\mu}|x-y|)&\mu G(\sqrt{\mu}|x-y|)
\end{pmatrix}
\end{equation*}
Therefore 
$$
|\rho_\mu(|x-y|)| =  \frac{1}{8|\mu|} \sqrt{2 |\mu|^2|G(\sqrt{\mu}|x-y|)|^2 + |\partial_{\bar{z}}^2 G(\sqrt{\mu}|x-y|)|^2 + |\partial_{z}^2 G(\sqrt{\mu}|x-y|)|^2 }.
$$
As a consequence, if we denote by $\rho_\theta(|x-y|)$ the kernel of the operator $(D_0-e^{i\theta})^{-1}$ then
\begin{equation}\label{rho1}
|\rho_\theta(r)|\leq C\ln r^{-1},\qquad {\rm if}\qquad r<1/2,
\end{equation}
and
\begin{equation}\label{rho2}
|\rho_\theta(r)|\leq Cr^{-1/2},\qquad {\rm if}\qquad r>1/2.
\end{equation}
In order  to prove  the latter  relations, one  has  to differentiate  the integral kernel of $(\Delta^2-\mu^2)^{-1}$, using the   formulas
$$
\frac{\partial r}{\partial z}=\frac12 \frac {\bar z}{r},\qquad  \frac{\partial^2 r}{\partial z^2}=-\frac14 \frac {\bar z^2}{r^3}
$$
and
$$
\frac{\partial r}{\partial \bar z}=\frac12 \frac { z}{r},\qquad  \frac{\partial^2 r}{\partial \bar z^2}=-\frac14 \frac { z^2}{r^3}.
$$
Since the integral kernel of $(\Delta^2-\mu^2)^{-1}$  is  $\frac{i}{8\mu} G(\sqrt\mu r)$, we  obtain  from \eqref{resolventm=0} that
\begin{multline*}
8|\rho_\theta(r)|\leq \Bigl( \Bigl| \frac{\partial^2G(e^{i\theta/2} r)}{\partial z^2} \Bigr|^2+   \Bigl| \frac{\partial^2G(e^{i\theta/2} r)}{\partial \bar z^2 }  \Bigr|^2 +2|G(e^{i\theta/2} r)|^2 \Bigr)^{1/2}\\
\leq  C(r^{-1}|G'(e^{i\theta/2} r)|+|G''(e^{i\theta/2} r)|+|G(e^{i\theta/2} r)|).
\end{multline*}
The  positive constants in the inequalities \eqref{rho1} and \eqref{rho2}  do not depend on $\theta\in[0,\pi/2]$. In particular,
$$
M:=\sup_\theta\sup_{x\in{\mathbb R}^2}\int_{{\mathbb R}^2}|\rho_\theta(|x-y|)|^qdy<\infty,\qquad q>4.
$$
Let us estimate now  the  norm of the operator
$T=W(D_0-e^{i\theta })^{-1}W$  with  the kernel
$$
\tau(x,y)=W(x)\rho_\theta(|x-y|)W(y).
$$
For that purpose, we  estimate the  sesquie-linear form
of this operator :
$$
(Tu,v)=\int_{{\mathbb R}^2}\int_{{\mathbb R}^2} \bar v(x)W(x)\rho_\theta(|x-y|)W(y)u(y)\,dxdy.
$$
Obviously,
$$
|(Tu,v)|^2=\Bigl|\int_{{\mathbb R}^2}\int_{{\mathbb R}^2} \bar v(x)W(x)\rho_\theta(|x-y|)W(y)u(y)\,dxdy\Bigr|^2\leq
$$
$$
\int_{{\mathbb R}^2}\int_{{\mathbb R}^2}  |v(x)|^2|\rho_\theta(|x-y|)||W(y)|^2\,dxdy\int_{{\mathbb R}^2}\int_{{\mathbb R}^2} |W(x)|^2|\rho_\theta(|x-y|)||u(y)|^2\,dxdy\leq
$$
$$
\Bigl(\sup_x\int_{{\mathbb R}^2} |\rho_\theta(|x-y|)||W(y)|^2\,dy\Bigr)^2\|u\|^2\,||v||^2\leq
$$
$$
\Bigl(\int_{{\mathbb R}^2} |\rho_\theta(|x-y|)|^q\,dy\Bigr)^{2/q} \|V\|_p^{2}\,\,\|u\|^2\,\|v\|^2,\qquad \frac1p+\frac1q=1,\qquad q>4.
$$
Therefore,
$$
\|T\|\leq C\|V\|_p,\qquad 1<p<4/3.
$$
We  are now able to estimate  the  norm of the operator $T_k=W(D_0-k)^{-1}W$  for $k\notin \sigma (D_0)$.
Indeed,
$$
|(T_ku,v)|=\left|\int_{{\mathbb R}^2}\int_{{\mathbb R}^2} \bar v(x)W(x)\rho_\theta(\sqrt{|k|}|x-y|)W(y)u(y)\,dxdy\right|=
$$
$$
\frac1{|k|^2}\left|\int_{{\mathbb R}^2}\int_{{\mathbb R}^2} \bar v(x/\sqrt{|k|})W(x/\sqrt{|k|})\rho_\theta(|x-y|)W(y/\sqrt{|k|})u(y/\sqrt{|k|})\,dxdy\right|\leq
$$
$$
\frac{C}{|k|^2}\|V(\cdot/\sqrt{|k|})\|_p  \|u(\cdot/\sqrt{|k|})\|\,\,\|v(\cdot/\sqrt{|k|})\|=\frac{C\|V\|_p}{|k|^{(p-1)/p}}\|u\|\,\,\|v\|.
$$
Consequently,
$$
\|T_k\|\leq \frac{C\|V\|_p}{|k|^{(p-1)/p}}.
$$

Observe now that the kernel of the operator $(\Delta^2-\mu^2)^{-1}$ is the function $iG(\sqrt{\mu}|x-y|)/(8\mu)$.
The function $G(\sqrt{\mu}|x-y|)$ has  the same properties as $\rho_\theta(\sqrt{|\mu|}|x-y|)$. Moreover it is  bounded.  Therefore,  by mimicking the  above arguments, one  proves  that
\begin{equation}\label{Delta2}
\|W(\Delta^2-\mu^2)^{-1}W\|\leq \frac{C\|V\|_p}{|\mu|^{(2p-1)/p}},\qquad 1\leq p<4/3.
\end{equation}
This  leads  to the  estimate
$$
\|W(D_m-k)^{-1}W\|^p\leq \frac{C\int_{{\mathbb R}^2} |V(x)|^p dx}{|\mu|^{p-1}}\left(\sqrt{\left|\frac{k+m}{k-m}\right|}+\sqrt{\left|\frac{k-m}{k+m}\right|}+1\right)^p,\qquad \mu^2=k^2-m^2.
$$
Now  the  statement of  our  theorem follows  from    the fact  that if $k$  is an eigenvalue of $D=D_m+V$, then $\|W(D_m-k)^{-1}W\|\geq1$.
The proof is complete.

\bigskip
\noindent
In the picture below we describe the areas of possible location of  complex eigenvalues  depending on the value of  $C\int_{\Bbb R^2} |V(x)|^p\, dx$, where $m=1$ and $p=1.2$.

\bigskip
\qquad\qquad \includegraphics[scale=0.50]{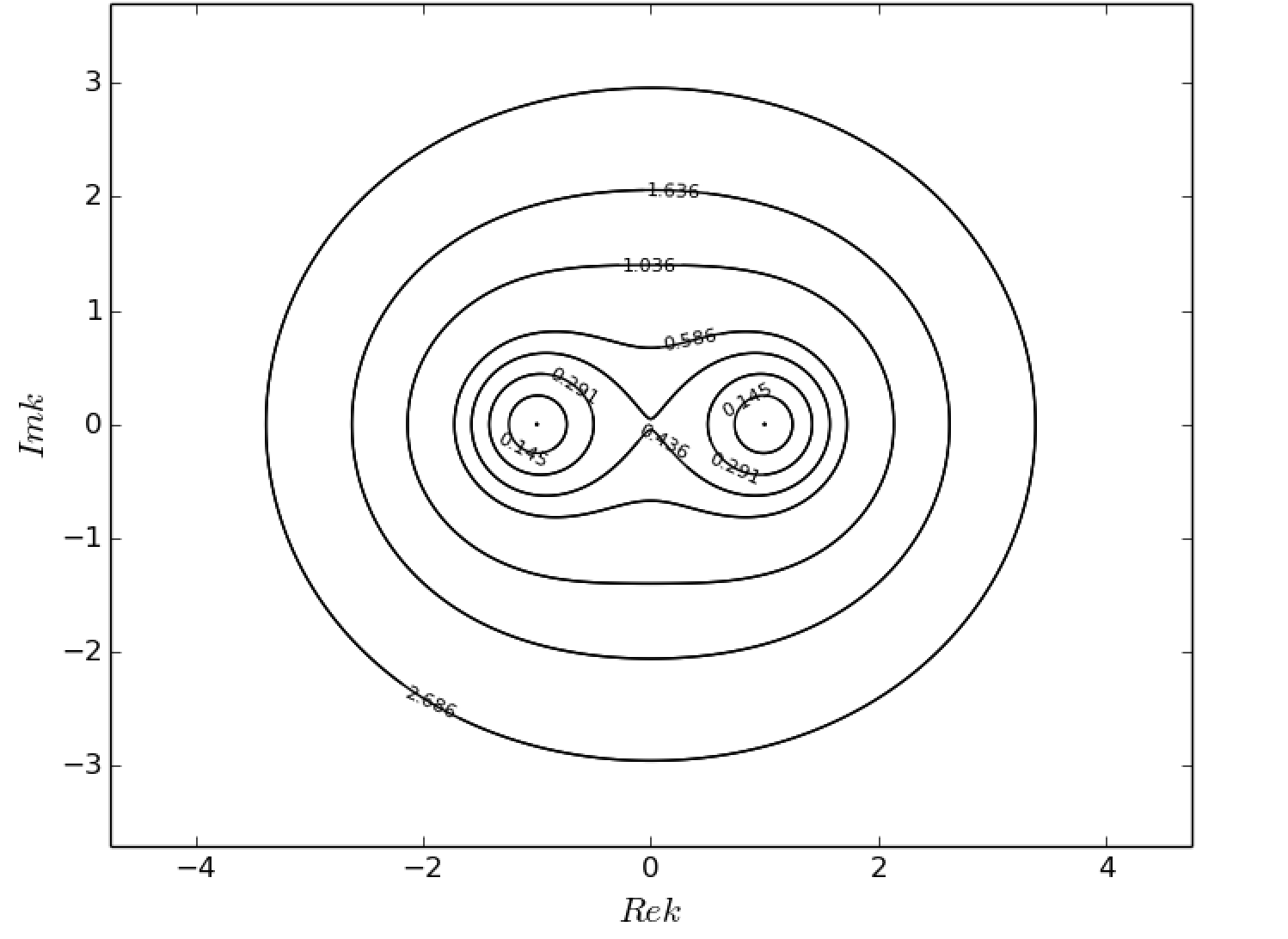}

\bigskip

\noindent
{\it Proof  of Theorem~\ref{t1.2}.}  As before we  use  the representation
$$
(D_m-k)^{-1}=(m\gamma_0+k-\mu)(\Delta^2-\mu^2)^{-1}+(D_0-\mu)^{-1}.
$$
The operator $(\Delta^2-\mu^2)^{-1}$  is an integral operator with the  kernel
$$
g_k(x,y)=\frac{i}{8\mu}\Bigl(H(\sqrt{\mu} r)-H(i\sqrt{\mu} r)\Bigr),
$$
where $H(z)=H^{(1)}_0(z)$  is  the Hankel function of first kind. 
 Again we denote
$$
G(z)=H(z)-H(iz).
$$
We need  to know the  behavior  of the  function $G$ only in the region $0<\arg z<\pi/2$, where
this function is  bounded.  The  boundedness of $G$  implies the estimate \eqref{Delta2}  with $p=1$.

It remains  to estimate   the  norm of the operator $T_\mu=W(D_0-\mu)^{-1}W$  for ${\rm Im }\, \mu> 0$.
We already know that if $\mu=|\mu|e^{i\theta}$ the operator $(D_0-\mu)^{-1}$ is an integral operator  with the kernel
$\rho_\theta(\sqrt{|\mu|}|x-y|)$, where $\rho_\theta$ is  a  function   having the properties 
\begin{equation}\label{r<z}
|\rho_\theta(r)|\leq C\ln r^{-1},\qquad {\rm if}\qquad r<1/2,
\end{equation}
and
\begin{equation}\label{r>z}
|\rho_\theta(r)|\leq Cr^{-1/2},\qquad {\rm if}\qquad r>1/2.
\end{equation}
The  positive constants in these inequalities  do not depend on $\theta\in[0,\pi]$.
As  before, we  estimate the  sesquie-linear form
of this operator :
$$
(T_\mu u,v)=\int_{{\mathbb R}^2}\int_{{\mathbb R}^2} \bar v(x)W(x)\rho_\theta(\sqrt{|\mu|}|x-y|)W(y)u(y)\,dxdy.
$$
Obviously,
$$
|(T_\mu u,v)|^2=\left|\int_{{\mathbb R}^2}\int_{{\mathbb R}^2} \bar v(x)W(x)\rho_\theta(\sqrt{|\mu|}|x-y|)W(y)u(y)\,dxdy\right|^2\leq
$$
$$
\int_{{\mathbb R}^2}\int_{{\mathbb R}^2}  |v(x)|^2|\rho_\theta(\sqrt{|\mu|}|x-y|)||W(y)|^2\,dxdy\int_{{\mathbb R}^2}\int_{{\mathbb R}^2} |W(x)|^2|\rho_\theta(\sqrt{|\mu|}|x-y|)||u(y)|^2\,dxdy\leq
$$
$$
\left(\sup_x\int_{{\mathbb R}^2} |\rho_\theta(\sqrt{|\mu|}|x-y|)||V(y)|\,dy\right)^2\|u\|^2\,\|v\|^2.
$$
Therefore,
$$
\|T_\mu \|\leq \sup_x\int_{{\mathbb R}^2} |\rho_\theta(\sqrt{|\mu|}|x-y|)||V(y)|\,dy.
$$
The bounds \eqref{r<z} and \eqref{r>z} imply
\begin{multline*}
 \sup_x\int_{{\mathbb R}^2} |\rho_\theta(\sqrt{|\mu|}|x-y|)||V(y)|\,dy\\
 \leq C
\left(|\ln|\mu||\sup_{x\in{\mathbb R}^2}\int_{|x-y|<(2|\mu|)^{-1}} |V(y)|\,dy+\sup_{x\in{\mathbb R}^2}\int_{{\mathbb R}^2}\left(1+|\ln|x-y||\right)|V(y)|dy\right),
\end{multline*}
which leads to
$$
\|T_\mu\|\leq C\left(|\ln|\mu||\sup_{x\in{\mathbb R}^2}\int_{|x-y|<(2|\mu|)^{-1}} |V(y)|\,dy+\sup_{x\in{\mathbb R}^2}\int_{{\mathbb R}^2}\left(1+|\ln|x-y||\right)|V(y)|dy\right).
$$
Since
$$
\|W(D_m-k)^{-1}W\|\leq 
\|W(m\gamma_0+k-\mu)(\Delta^2-\mu^2)^{-1}W\|+\|T_\mu\|
$$
and since \eqref{Delta2} holds with $p=1$, we  obtain
\begin{multline*}
\|W(D_m-k)^{-1}W\|\\
\leq C\left(|\ln|\mu||\sup_{x\in{\mathbb R}^2}\int_{|x-y|<(2|\mu|)^{-1}} |V(y)|\,dy+\sup_{x\in{\mathbb R}^2}\int_{{\mathbb R}^2}\left(1+|\ln|x-y||\right)|V(y)|dy
\right)\\
+C\int_{{\mathbb R}^2}|V(x)|dx\left(\sqrt{\left|\frac{kt-m}{k+m}\right|}+\sqrt{\left|\frac{k+m}{k-m}\right|}+1\right).
\end{multline*}
The  statement of   Theorem~\ref{t1.2} follows  from    the fact  that if $k$  is an eigenvalue of $D=D_m+V$, then $\|W(D_m-k)^{-1}W\|\geq1$.

\bigskip

\section{A special case}

Consider  now  a special case, when $V(x)=iW^2(x)$, where $W(x)=W^*(x)$ is  a matrix  valued  function.  It  turns out,  that in this case we can get a more precise information about the spectral properties of the operator $D$. 
\begin{theorem}
Let $k\notin \sigma(D_m)$ be an eigenvalue of the operator $D=D_m+V$, where $V=iW^2$.  Let $\mu$  be the number in  the upper half-plane  defined  by $\mu^2=k^2-m^2$.
Then
\begin{equation}\label{D+iV}
\left(C\left(\left| \frac{k+m}{\mu}-1\right|+\left| \frac{k-m}{\mu}-1\right|\right)+1\right)\, \frac14\int_{{\Bbb R}^2}{\rm tr}|V|dx\geq1,
\end{equation}
where  the constant $C$ is independent of $V$, $m$ and $k$.
\end{theorem}

\bigskip

{\it Proof.}
According to the Birman-Schwinger principle, $k$  is an eigenvalue  of the operator $D$ if and  only if  $1$ is an eigenvalue of the operator $X=-iW(D_m-k)^{-1}W$.
On the other  hand, if  $1$ is an eigenvalue of $X$  then $\|{\rm Re }X\|\geq 1$.
Since, 
$$
{\rm Re }X=W\, {\rm  Im}(D_m-k)^{-1}W,
$$
we would like to have the explicit expression for the operator
$
{\rm Im }(D_m-k)^{-1}.
$
Let us first obtain this  representation for the  case $m=0$. Observe  that  $D_0$ is the operator  with the symbol   
$$
\begin{pmatrix}
0 & -(\xi_1+i\xi_2)^2\\ -(\xi_1+i\xi_2)^2&0
\end{pmatrix}.
$$
The eigenvalues of this matrix are $\pm |\xi_1+i\xi_2|^2$. The orthogonal  projections $P_1(\xi)$ and $ P_2(\xi)$, $\xi=\xi_1+i\xi_2$,  onto  the  eigenvectors  depend  only on $\arg (\xi)$. Therefore,  the symbol of $D_0$ is
$$
|\xi|^2P_1(\xi)-|\xi|^2P_2(\xi),
$$
which implies  that the integral kernel of ${\rm Im}\,(D_0-k)^{-1}$ is
$$
(2\pi)^{-2}\int_{{\Bbb R}^2}\exp(i\xi(x-y))\left(\frac{({\rm Im} \,k)\, P_1(\xi)}{(|\xi|^2-{\rm Re} \, k )^2+({\rm Im}\, k)^2}+\frac{ ({\rm Im} \,k)\,  P_2(\xi)}{(-|\xi|^2- {\rm Re}\, k)^2+({\rm Im}\, k)^2}\right)\,\,d\xi.
$$
It  follows  from this  representation that 
the kernel  of the operator ${\rm Im}\,(D_0-k)^{-1}$ is  bounded  by $1/4$ as using polar coordinates and changing variables $|\xi|^2 = t$ we obtain
$$
\|{\rm Im}\,(D_0-k)^{-1}\| \le (2\pi)^{-2}\, \int_{\Bbb S^1} \int_{-\infty}^\infty \frac {{\rm Im} \,k}{t^2 + ({\rm Im} \,k)^2} \, dt = \frac14, \qquad {\rm Im} \,k>0.
$$
Consequently,
$$
\left\|W{\rm Im}\,(D_0-k)^{-1}W\right\| \leq   {\rm tr}\, \left(W{\rm Im}\,(D_0-k)^{-1}W \right)\leq \frac14\int_{{\Bbb R^2}}{\rm tr}\,W^2(x)dx.
$$
If $m>0$, then we have
$$
\| {\rm Re}X\| \leq \left\|\frac1{2\mu}W(m\gamma_0+k-\mu)(\Delta^2-\mu^2)^{-1}W \right\|+\left\| W\, {\rm Im}(D_0-\mu)^{-1}W\right\|
$$
and  that
$$
\left\|\frac{m\gamma_0+k-\mu}{2\mu}\right\|\leq \frac12\left(\left| \frac{k+m}{\mu}-1\right|+\left| \frac{k-m}{\mu}-1\right|\right).
$$
It  remains to  note that according to \eqref{Delta2} with $p=1$,
$$
\|W(\Delta^2-\mu^2)^{-1}W\|\leq C \int_{{\Bbb R}^2}{\rm tr |V|}\,dx
$$
The  proof is completed. $\,\,\,\,\,\,\,\,\,\,\,\,\,\,\, \Box$

\bigskip

The next  result  says  that   the spectrum of the operator $D_0$  is  stable with  respect  to
  small perturbations  of the   form $V=iW^2$. 

\begin{corollary} Let $m=0$ and
let $V=iW^2$ with $W^*=W$.  
Assume  that
\begin{equation}\label{m=0}
\frac14\int_{{\Bbb R}^2}{\rm tr}|V|dx<1.
\end{equation}
 Then   the operator $D=D_0+V$  does not have eigenvalues outside of the real line ${\Bbb R}$, i.e.  the spectrum of $D$ is real.
\end{corollary}

\bibliographystyle{amsalpha}

\end{document}